%% file: Landaufinal.tex
\documentclass[12pt]{amsart}
\usepackage{amsmath,amsfonts,amsthm,amsopn,cite,mathrsfs}
\usepackage{epsfig,verbatim}
\include{macros}

\newcommand{\Ghat}{\widehat{G}}
\newcommand{\bom}{\mathcal{B}_\Omega}
\newcommand{\Z}{{\Bbb Z}}
\begin{document}
\begin{abstract}
We derive necessary conditions for sampling and interpolation of
bandlimited functions on a locally compact abelian group in line
with the classical results of H. Landau for bandlimited functions
on ${\Bbb R}^d$. Our conditions are phrased as comparison
principles involving a certain canonical lattice.
\end{abstract}

\title[Landau's necessary density conditions]{Landau's necessary density
conditions \\ for LCA groups}
\author{Karlheinz Gr\"ochenig}
\address{Faculty of Mathematics \\
University of Vienna \\
Nordbergstrasse 15 \\
A-1090 Vienna, Austria}
\email{karlheinz.groechenig@univie.ac.at}
\thanks{K.G. was supported by the  Marie-Curie Excellence Grant
  MEXT-CT-2004-517154}

\author{Gitta Kutyniok}
\address{Department of Statistics\\
Stanford University\\
Stanford, CA 94305, USA} \email{kutyniok@stanford.edu}
\thanks{G.K. was supported by Preis der Justus-Liebig-University Giessen 2006
and Deutsche Forschungsgemeinschaft (DFG) Heisenberg fellowship KU
1446/8-1.}

\author{Kristian Seip}
\address{Department of Mathematical Sciences\\
Norwegian University of Science and Technolgy\\
NO-7491 Trondheim, Norway}
\email{seip@math.ntnu.no}
\thanks{K.S. was supported by the Research Council of Norway grant 10323200.}
\thanks{This research is part of the European Science Foundation Networking
Programme Harmonic and Complex Analysis and Its Applications
(HCAA)} 
\date{\today}

\maketitle

\section{Introduction}

H. Landau's necessary density conditions for sampling and
interpolation \cite{Lan67} may be viewed as a general principle
resting on a basic fact of Fourier analysis: The complex
exponentials $e^{i kx}$ ($k$ in $\mathbb{Z}$) constitute an
orthogonal basis for $L^2([-\pi,\pi])$. The present paper extends
Landau's conditions to the setting of locally compact abelian
(LCA) groups, relying in an analogous way on the basics of Fourier
analysis. The technicalities---in either case of an operator
theoretic nature---are however quite different. We will base our
proofs on the comparison principle of J. Ramanathan and T. Steger.
\cite{RS95}.

We recall briefly Landau's results, suitably adapted to our
approach. Let $\Omega$ be a bounded measurable set in ${\Bbb R}^d$
and let $\bom$ denote the subspace of $L^2({\Bbb R}^d)$ consisting
of those functions whose Fourier transform is supported on
$\Omega$. We say that a subset $\Lambda$ of $\bR^d$ is
\emph{uniformly discrete} if the distance between any two points
exceeds some positive number. A uniformly discrete set $\Lambda$
is a \emph{set of sampling for} $\bom$ if there exists a constant
$C$ such that $\|f\|^2_2\le C\sum_{\lambda\in \Lambda}
|f(\lambda)|^2$ for every $f$ in $\bom$, and a \emph{set of
interpolation for} $\bom$ if, for each square-summable sequence
$\{a_\lambda\}_{\lambda\in \Lambda}$, there is a solution $f$ in
$\bom$ to the interpolation problem $f(\lambda)=a_\lambda$,
$\lambda$ in $\Lambda$.

The canonical case is when $\Omega$ is a cube of side length
$2\pi$ and $\Lambda$ the integer lattice $\Z^d$. Since the complex
exponentials $e^{i \lambda\cdot x}$ ($\lambda$ in $\Lambda$)
constitute an orthogonal basis for $L^2(\Omega)$, it is immediate
by the Plancherel identity that $\Lambda$ is both a set of
sampling and a set of interpolation for $\bom$. This result scales
in a trivial way: $c\Z^d$ is a set of sampling and a set of
interpolation for $\bom$ when $\Omega$ is a cube of side length
$c^{-1}2\pi$. If we agree that the density of the integer lattice
is $1$, then we have that the density of the lattice equals
$(2\pi)^{-d}$ times the volume of the spectrum $\Omega$.

Landau's work may be understood as saying that this density result
takes the following form for general $\Omega$ and uniformly
discrete sets $\Lambda$: \emph{\begin{itemize} \item[(S)] If
$\Lambda$ is a set of sampling for $\bom$, then $\Lambda$ is
everywhere at least as dense as the lattice
$(2\pi)^{-1}|\Omega|^{1/d}\, \Z^d$. \item [(I)]If $\Lambda$ is a
set of interpolation for $\bom$, then $\Lambda$ is everywhere at
least as sparse as the lattice $(2\pi)^{-1} |\Omega|^{1/d}\,
\Z^d$.
\end{itemize}}

Landau gave precise versions of these statements in terms of the
following notion of density. For $h
> 0$ and $x$ a point in $\bR^d$, let $Q_h(x)$ denote the closed cube
centered at $x$ of side length $h$. Then the \emph{lower Beurling
density} of the uniformly discrete set $\Lambda$ is defined as
\[\cD^-_B(\La) = \liminf_{h \to \infty}\inf_{x \in \bR^d}
    \dfrac{ \mathrm{card} \, (\La \cap Q_h(x))} {h^d} \, ,\]
and its \emph{upper Beurling density} is
\[\cD^+_B(\La)= \limsup_{h \to \infty}\sup_{x \in \bR^d}
    \dfrac{\mathrm{card} \, (\La \cap Q_h(x))} {h^d}.\]
Landau's result says that a set of sampling $\Lambda$ for $\bom$
satisfies $\cD^-_B(\La)\ge (2\pi)^{-d}|\Omega|$ and a set of
interpolation $\Lambda$ for $\bom$ satisfies $\cD^+_B(\La)\le
(2\pi)^{-d}|\Omega|$.

Given two uniformly discrete sets $\Lambda$ and $\Lambda'$ and
nonnegative numbers $\alpha$ and $\alpha'$, we write $\alpha
\Lambda\preceq \alpha' \Lambda'$ if for every positive $\epsilon$
there exists a compact subset $K$ of $\rd$ such that for every
compact subset $L$ we have \begin{equation} \label{dominateRd}
(1-\epsilon) \,\alpha\, \mathrm{card}\, (\Lambda \cap L )
  \leq  \alpha'\, \mathrm{card}\, (\Lambda' \cap (K+L)).
  \end{equation}
With this notation, we have the following equivalent way\footnote{A
proof of the equivalence is given in Section \ref{sec:dense} of this
paper.} of expressing Landau's density conditions:
\emph{\begin{itemize}\item[(S)] If $\Lambda$ is a set of sampling
for $\bom$, then $(2\pi)^{-d}\, |\Omega| \,\Z^d\, \preceq\, \La $.
 \item[(I)] If $\Lambda$ is a set of interpolation for $\bom$,
then $\La \, \preceq \,(2\pi)^{-d}\,|\Omega|\, \Z^d $.
\end{itemize}}
The latter formulation may look less appealing than that given by
Landau, but it has the advantage of presenting Landau's density
conditions as a \emph{comparison principle}; we note that this
version does not require the use of dilations of cubes, which in
general LCA groups make no sense.

We take as our starting point this reformulation of Landau's
results, and in brief our plan is as follows. We need to identify,
in the general setting of LCA groups, a canonical case to be used
for comparison. Then, besides comparing $\Lambda$ with a suitable
canonical lattice, we also need to compare spectra. We will do this
by estimating a general spectrum $\Omega$ in terms of a disjoint
union of small ``cubes''. A nontrivial point will be to clarify what
are the right ``cubes'' and what are the ``lattices'' associated
with such sets. The technicalities of the comparison will in fact
take place at this ``atomic'' level, and it is here that the
Ramanathan--Steger comparison lemma will play a crucial role.

While our approach leads to a best possible asymptotic result in
the more general setting of LCA groups, we lose a subtle level of
precision compared to Landau's work, which is based on estimates
for the eigenvalues of a certain concentration operator. For
$\Omega$ a finite union of real intervals, Landau obtained sharp
bounds for the number of points from a set of sampling or of
interpolation to be found in a large interval $I$. In these bounds
appears an additional term of order $\log|I|$, and---as shown in
\cite{OS02}---this can be seen as a manifestation of the
John--Nirenberg theorem.

It is worth noting that one may encounter situations in which no
obvious analogue of a lattice is available. An interesting example
is that of the unit sphere in ${\Bbb R}^d$. In a recent paper
\cite{Mar06}, J. Marzo managed to employ Landau's method in this
setting without any explicit comparison between uniformly discrete
sets. In our setting, the group of $p$-adic numbers is an example
of an LCA group that fails to contain a lattice. Our approach will
be to restrict to a discrete quotient on which a meaningful
comparison with a lattice can be made.

\medskip


The ideas of Ramanathan and Steger have been employed by many
authors. We would in particular like to mention the basic theory
developed by R. Balan, P. Casazza, C. Heil, and Z. Landau in
\cite{BCHL06}. That paper introduces a notion of density for
frames parameterized by discrete abelian groups, such as Gabor
frames. The present paper is however only loosely related to
\cite{BCHL06}; we will require a more general notion of density,
since we will be dealing with uniformly discrete sets in general
LCA groups rather than discrete abelian groups.

\section{Landau's density theorem for LCA groups}
\label{sec:lcagroups}

We start by recalling some basic facts about locally compact abelian
(LCA) groups. For more information we refer to the books
\cite{Fol95} and \cite{HR63}.

Let $G$ be a locally compact abelian group; to avoid trivialities,
we assume that $G$ is \emph{non-compact}. The group multiplication
will be written multiplicatively as $xy$, and we will use the
notation $xK = \{ y\in G: y= xk, k\in K\}$ and $KL = \{ y\in G: y=
kl, k\in K, l\in L\}$ for $x \in G$ and $K$, $L$ being subsets of
$G$. A locally compact group $G$ is always equipped with a Haar
measure, which in the following will be denoted by $\mu_G$. We
follow the convention that the Haar measure of a compact (sub)group
is normalized to be a probability measure.


Let $\widehat{G}$ be the dual group of $G$. We write the action of
a character $\omega \in \Ghat $ on $x\in G$ by $\langle \omega
,x\rangle$. The {\em annihilator} $H^\perp \subseteq \Ghat $ of a
subgroup $H\subseteq G$ is defined as $H^\perp = \{ \chi \in \Ghat
: H \subseteq \mathrm{ker}\, \chi \}$. By Pontrjagin duality, we
can identify $G $ with $\widehat{\Ghat}$, and we will frequently
use that $\widehat{H} \simeq \Ghat / H^\perp $ and $(G/H)\, \widehat{}
\, \simeq H^\perp $.

The Fourier transform $\cF $ is  defined by
\[
  \cF f(\omega)=\widehat{f}(\omega ) = \int _G f(x) \overline{\langle\omega ,x\rangle}
  \, d\mu_G(x),  \qquad
  \omega \in \Ghat \, .
\]
We assume that the Haar measures on $G$ and $\Ghat $ have been
chosen such that $\cF $ is a unitary map from $L^2(G)$ onto
$L^2(\Ghat )$, in accordance with Plancherel's theorem. If $\Omega
\subseteq \Ghat $ is a measurable set of positive measure,
\[
  \bom = \{ f\in L^2(G) : \supp \, \widehat{f} \subseteq \Omega \}
\]
is the space of ``band-limited'' functions with spectrum in
$\Omega $.

A subset $\Lambda$ of $G $ is called {\em uniformly discrete} if
there exists an open set $U$ such that the sets $\lambda U$
($\lambda$ in $\Lambda$) are pairwise disjoint. The definition of
sets of sampling and interpolation given in the introduction extends
without any change to the setting of LCA groups. We are interested
in such sets for the space $\bom$.

We will assume that the dual group $\widehat{G}$ is compactly
generated. This may seem a rather severe restriction and means
that for instance $p$-adic groups are excluded from our
consideration. However, if the spectrum is relatively compact, we
may assume without loss of generality that $\Ghat$ is compactly
generated. For a clarification of this point, we refer to
Section~\ref{sec:red}. By the structure theory of LCA
groups~\cite{HR63}, $\Ghat $ is then isomorphic to $\rd \times \bZ
^m \times K_0$ for  a compact group $K_0$. Consequently,  $G$ is
of the form $G= \rd \times \bT
 ^m \times D_0$ with $D_0$ a (countable)  discrete  group. We then select the
 uniformly discrete subset $\Gamma_0=\bZ^d\times\{e\}\times D_0$ as
 the canonical lattice to be used for comparison, where $e$ is the identity
 element in $\bT^m$. We assume that the Haar measure $\mu_{\Ghat}$
 is normalized so that $\mu_{\Ghat} ([-\pi,\pi]^d \times \{e\} \times K_0)=1$.

We define the relation `$\preceq$' for uniformly discrete subsets of
$G$ as we did in \eqref{dominateRd}: Given two uniformly discrete
sets $\Lambda$ and $\Lambda'$ and nonnegative numbers $\alpha$ and
$\alpha'$, we write $\alpha \Lambda\preceq \alpha' \Lambda'$ if for
every positive $\epsilon$ there exists a compact subset $K$ of $G$
such that for every compact subset $L$ we have \begin{equation}
\label{dominatelca} (1-\epsilon) \,\alpha\,  \mathrm{card}\,
(\Lambda \cap L )
  \leq  \alpha'\, \mathrm{card}\, (\Lambda' \cap KL). \end{equation}
It is immediate from the definition given by \eqref{dominatelca}
that the relation `$\preceq$' is transitive, a fact that will be
used repeatedly in what follows.

With this notation, we may state Landau's necessary conditions for
sampling and interpolation in the context of a general LCA group
as follows.

\begin{tm} \label{tm:compc}
Suppose $\Lambda$ is a uniformly discrete subset of the LCA group
$G$ and $\Omega$ is a relatively compact subset of $\widehat{G}$.
\begin{itemize}
\item[(S)] If $\Lambda $ is a set of sampling for $\bom $, then
$ \mu_{\Ghat}(\Omega)\, \Gamma_0\, \preceq \,
 \Lambda$. \item[(I)] If $\Lambda $ is a set of interpolation
for $\bom $, then  $  \Lambda\,  \preceq \,
 \mu_{\Ghat}(\Omega)\, \Gamma_0$.
\end{itemize}
\end{tm}

One may think of $\mu_{\Ghat}(\Omega)$ as the ``Nyquist density''.
Indeed, the relation `$\preceq$' gives us a way of defining
densities of a uniformly discrete set: The \emph{lower uniform
density} of $\La$ is defined as \[ \cD^-(\La)=\sup\{\alpha: \
\alpha \Gamma_0 \preceq \Lambda\}, \] and its \emph{upper uniform
density} is
\[ \cD^+(\La)=\inf\{\alpha: \ \Lambda\preceq \alpha \Gamma_0\},
\]
with the understanding that $\cD^+(\La)=\infty$ if the set on the
right hand side is empty. We will later show that both densities are
always finite, and so the infimum in the definition of
$\cD^{-}(\Lambda)$ is in fact a minimum, and the supremum in the
definition of $\cD^{+}(\Lambda)$ is a maximum. With these
definitions, Theorem~\ref{tm:compc} can be reformulated in the
following classical way.
\begin{Th1prime} 
  Suppose $\Lambda$ is a uniformly discrete subset of the LCA group
$G$ and $\Omega$ is a relatively compact subset of $\widehat{G}$.
\begin{itemize}
\item[(S)] If $\Lambda $ is a set of sampling for $\cB _\Omega$,
then
  $\mathcal{D} ^-(\Lambda ) \geq \mu _{\widehat{G}}(\Omega )$.
  \item[(I)] If
  $\Lambda $ is a set of interpolation for $\cB _\Omega$, then
  $\mathcal{D} ^+(\Lambda ) \leq \mu _{\widehat{G}}(\Omega )$.
\end{itemize}
\end{Th1prime}
We will show below (Lemma~\ref{equiv}) that when $G=\rd$, $\cD^-(\Lambda)$ and
$\cD^+(\Lambda)$ reduce to the usual Beurling densities. Indeed,
we will see that an ``intermediate'' formulation of the densities,
valid for any LCA group $G$, may be obtained by replacing the
counting measure of $\Gamma_0$ by the Haar measure $\mu_G$. We
will also show that, in general, $\cD^-(\Lambda)\le
\cD^+(\Lambda)<\infty$. A particular consequence of this bound is
that $\cD^-(\Gamma_0)=\cD^+(\Gamma_0)=1$, because the transitivity
of the relation `$\preceq$' implies that either
$\cD^-(\Gamma_0)=\cD^+(\Gamma_0)=1$ or
$\cD^-(\Gamma_0)=\cD^+(\Gamma_0)=\infty$.

We will return to this discussion of uniform densities in
Section~\ref{sec:dense}, after the proof of Theorem~\ref{tm:compc}.
That proof requires some preparation, to be presented in the next
three paragraphs. The most significant ingredients are the Fourier
bases for small ``cubes'', given in Section~\ref{sec:cubes}, and the
Ramanathan--Steger comparison principle, treated in
Section~\ref{sec:comparison}. The actual proof of
Theorem~\ref{tm:compc} is given in Section~\ref{sec:proof}.

After a consideration of the case when $\Ghat$ is not compactly
generated in Section~\ref{sec:red}, we close in
Section~\ref{sec:fin} with some additional remarks pertaining to
Theorem~\ref{tm:compc}.

\section{Square sums of point evaluations at uniformly discrete sets}
\label{sec:uds}

The purpose of this section is mainly to show that our a priori
assumption that $\Lambda$ be a uniformly discrete set implies no
loss of generality. However, one piece of this discussion will be
needed in the proof of Theorem~\ref{tm:compc}. This is
Lemma~\ref{help} below, which says that uniformly discrete sets
generate Carleson measures in a natural way.

We may of course remove the a priori assumption that a set of
interpolation be uniformly discrete, but it is easy to see that,
at any rate, a set of interpolation will be uniformly discrete.
The argument is standard. We first note that we can always solve
the interpolation problem with control of norms. This means that
if $\Lambda$ is a set of interpolation, then there exists a
constant $M$ such that the interpolation problem
$f(\lambda)=a_\lambda$ can be solved with $f$ in $\bom$ in such a
way that \[ \|f\|_2^2\, \le \, M \,\sum_\lambda |a_\lambda|^2\]
for every square-summable sequence
$\{a_\lambda\}_{\lambda\in\Lambda}$. This well-known fact is a
consequence of the open mapping theorem. Now assume that for every
open set $U$ in $G$ there are points $\lambda_1$ and $\lambda_2$
in $G$ such that $\lambda_1^{-1}\lambda_2 $ is in $U$. Solving the
problem $f(\lambda_1)=1$ and $f(\lambda)=0$ for every other
$\lambda$ in $\Lambda$, we get that $\|f\|_2\leq M$ and
\begin{eqnarray*}
  1 &=& |f(\lambda _1) - f(\lambda _2)| \leq \int _\Omega
 |\hat{f}(\omega )| |\langle \omega , \lambda _1\rangle - \langle
 \omega , \lambda _2\rangle| d\mu _{\hat{G}}(\omega)\\
&  \le & M \mu
 _{\hat{G}} (\Omega )^{1/2} \sup_{\omega\in \Omega}
|1-\langle \omega, \lambda^{-1}_1\lambda_2\rangle|,
\end{eqnarray*}
 which cannot
hold for arbitrary $U$ when $\Omega$ is relatively compact.

The reduction from a more general definition of sets of sampling
follows the same pattern as in \cite[pp. 140--141]{S95}. We will
therefore be brief and only mention a few technical modifications.
We begin with the following result on Carleson measures.

\begin{lemma}\label{help}
Let $\Lambda$ be a uniformly discrete subset of $G$, and assume
$\Omega$ is a relatively compact subset of $G$.  Then there is a
positive constant $C$ such that
\[\sum_{\lambda \in
\Lambda} |f(\lambda )|^2 \le C \|f\|_2^2
\]
holds for every $f$ in $\bom$.
\end{lemma}

\begin{proof}
  The proof is identical to that of Lemma~1 in~\cite{GR96}. Choose
  a function $g$ in $ L^1(G)$  so that $\widehat{g} (\omega )  = 1$ for
  $\omega \in \overline{\Omega}$ and such that for any (symmetric) compact
  neighborhood $U$ of $e$, the function $g^\sharp(x) = \sup _{u\in U }
  |g(xu)|$ is also in $L^1(G) $. (Such a function exists by
  ~\cite{Rei68}.) If $f$ is in  $\bom$, then $f=f\ast g$ and
  $f^\sharp (x) \leq (|f| \ast g^\sharp )(x) $ for all $x$ in
  $G$. Consequently, $\|f^\sharp \|_2 \leq \|f\| _2 \| g^\sharp \|_1  $
  for all $f$ in $\bom $. Clearly, $|f(\lambda )| \leq f^\sharp (x)$
  whenever $x \in \lambda U$. Since $\Lambda$ is uniformly discrete, we may choose
  $U$ such that
  \begin{eqnarray}
    \sum _{\lambda \in \Lambda}|f(\lambda )|^2
    &=& \sum _{\lambda \in \Lambda} \frac{1}{\mu_G(U)} \int _{\lambda
        U} |f(\lambda )|^2 d\mu _G(x) \notag  \\
&\leq & \sum _{\lambda \in \Lambda} \frac{1}{\mu_G(U)} \int
_{\lambda
        U} |f^\sharp (x)|^2 d\mu _G(x) \label{eq:cd01}\\
&\leq & \frac{1}{\mu_G(U)} \int _{G} |f^\sharp (x)|^2 d\mu
_G(x)\le \frac{\|g^\sharp\|_1^2}{\mu_G(U)} \|f\|_2^2. \notag
  \end{eqnarray}
  \end{proof}
This lemma and the $G$-invariance of $\bom$ imply that an
inequality of the form
\[ \sum_{\lambda \in
\Lambda} |f(\lambda )|^2 \le C \|f\|_2^2, \] valid for every $f$
in $\bom$, holds if and only if $\Lambda$ is a finite union of
uniformly discrete sets. The existence of such an inequality is
sometimes explicitly required in the definition of a set of
sampling.

We may now go one step further and prove that if there are positive
constants $c$ and $C$ such that
\[  c \|f\|_2^2 \le \sum_{\lambda \in
\Lambda} |f(\lambda )|^2 \le C \|f\|_2^2 \] holds for every $f$ in
$\bom$, then there are a uniformly discrete subset $\Lambda'$ of
$\Lambda$ and positive constants $c'$ and $C'$ such that
\[  c' \|f\|_2^2 \le \sum_{\lambda' \in
\Lambda'} |f(\lambda' )|^2 \le C' \|f\|_2^2 \] for every $f$ in
$\bom$. The key ingredient in the proof of this result is the
following continuity property: Suppose $\Lambda$ is a uniformly
discrete subset of $G$. Then, for every $\varepsilon>0$, there
exists a neighborhood $U$ of the identity $e$ such that if
$\lambda\mapsto \lambda'$ is a mapping from $\Lambda$ to $G$
satisfying $\lambda'\lambda^{-1}\in U$, then we have
\begin{equation} \label{bernstein}  \sum_{\lambda \in \Lambda}
|f(\lambda)-f(\lambda')|^2\le \varepsilon \|f\|_2^2
\end{equation} for every $f$ in $\bom$.

We give the short proof of \eqref{bernstein} and refer otherwise
to Lemma 3.11 of \cite{S95}. We let $g$ be as in the proof of
Lemma~\ref{help} and obtain
\[
    \sum _{\lambda \in \Lambda}|f(\lambda )-f(\lambda')|^2
    \le  \sum_{\lambda \in \Lambda} \left(\int_G |f(y)||g(\lambda y^{-1})
    -g(\lambda'y^{-1})|d\mu_G(y) \right)^2 \ \ \ \ \ \ \ \ \ \ \ \ \ \ \]
\[ \ \ \ \ \ \ \ \  \ \ \ \ \le \sum_{\lambda \in \Lambda} \int_G |f(y)|^2|g(\lambda
y^{-1})
    -g(\lambda'y^{-1})|d\mu_G(y)
    \int_G |g(\lambda x^{-1})-g(\lambda'x^{-1})|d\mu_G(x).\]
Since the translation operator $g(x)\mapsto g(\xi x)$ is continuous with respect to
the $L^1$-norm, the integral to the right can be made arbitrarily
small by a suitable choice of $U$, which is an estimate that is
uniform with respect to $\lambda$ and $\lambda'$. In the integral
to the left, we may then interchange the order of summation and
integration and essentially repeat the calculation made in
\eqref{eq:cd01} with $g$ in place of $f$. With a suitable choice
of $U$, the resulting estimate is \eqref{bernstein}.

\section{Fourier bases on small ``cubes''}\label{sec:cubes}

We will in what follows rewrite sampling and interpolation
properties in terms of the spanning properties of the resulting
functions on the Fourier transform side. By Lemma~\ref{help}, if
$\Omega$ is relatively compact, then $\Lambda$ is a set of
sampling for $\bom$ \fif\ the system $\{ e_\lambda (\omega ) =
\langle \omega ,\lambda \rangle \chi _\Omega (\omega ): \lambda
\in \Lambda \}$ is a frame for $L^2(\Omega)\subseteq L^2(\Ghat )$.
Likewise, $\Lambda $ is a set of interpolation for $\bom$ \fif\
$\{ e_\lambda (\omega ) = \langle \omega , \lambda \rangle \chi
_\Omega (\omega ): \lambda \in \Lambda \}$ is a Riesz sequence in
$L^2 (\Omega )$. This means that $\{ e_\lambda :\lambda \in
\Lambda \}$  is a Riesz basis in the closed linear span of the
functions $\{e_\lambda \}$.

We may at once apply this observation to the canonical lattice
$\Gamma_0=\Z^d\times\{e\}\times D_0$ of Theorem~\ref{tm:compc}.
Indeed, writing as before $\Ghat=\rd\times \Z^m \times K_0$, we
note that the characters labelled by $\Gamma_0$ and restricted to
$\Omega_0:=[-\pi,\pi]^d\times\{e\}\times K_0$ constitute an
orthonormal basis for $L^2(\Omega_0)$. Consequently, $\Gamma_0$ is
both a set sampling and a set of interpolation for
$\cB_{\Omega_0}$. (See also~\cite{kluvanek}.)

In the classical case when $G=\rd$, this is all we need, because
we can just scale $\Gamma_0$ to obtain Fourier bases for
arbitrarily small cubes\footnote{We recall from the introduction
that the motivation for such a rescaling is that we wish to
approximate an arbitrary spectrum by a union of small ``cubes''.}.
For general LCA groups, we need a different approach. It is
convenient to introduce some notation in order to state the lemma
to be used in place of a simple rescaling. We will say that a
discrete subgroup $\Gamma$ of $G$ is a {\em lattice} if the
quotient
$G/\Gamma$ is compact. 
A uniformly discrete set $\Gamma$ in $G$ will be said to be a
\emph{quasi-lattice} if the following holds. There is a compact
subgroup $K$ of $\widehat{G}$ and a lattice $\Upsilon$ in
$K^\perp$ such that $\Gamma=\{\widehat{k}\upsilon\}$, where
$\upsilon$ ranges over $\Upsilon$ and $\widehat{k} \in G$ ranges
over a set of representatives of $G/K^\perp $ in $G$. We may
identify $\{ \widehat{k} \}$ with $\widehat{K} \simeq G/K^\perp $,
and consequently $\{ \langle k, \widehat{k} \rangle\} $ ($k$ in
$K$) is an orthonormal basis for $L^2(K, \mu _K)$. We note that
every lattice $\Lambda$ is in particular a quasi-lattice; just
take $K=\{e\}$ and $\Upsilon=\Lambda$.




\begin{lemma}   \label{cubes}
Let $G$ be an LCA group whose dual group $\Ghat$ is compactly
generated. For every open neighborhood $U$ of the identity $e$ in
$\Ghat$ there exists a
  relatively compact subset $C$ of $U$  and a quasi-lattice
  $\Gamma$ in $G$
  with the following properties:
\begin{itemize}
\item[(i)] $L^2(C)$ possesses an orthogonal basis of characters
restricted to $C$ and labelled by $\Gamma$. \item[(ii)] There
exists a discrete subset $D$ of $\Ghat $ such that the translates
$dC, d\in D,$ form a partition of $\Ghat $.
\end{itemize}
\end{lemma}

\begin{proof}
   Since $\Ghat$
is compactly generated,  the structure theory implies that  any
neighborhood $U\subseteq \Ghat $ of $e$ contains a compact
subgroup $K$, such that $H := \Ghat /K \simeq  \rd \times \bZ ^m
\times \bT ^\ell \times F$, where $F$ is a  finite group and $d, m
, \ell \geq 0$. See~\cite[Thm.~9.6]{HR63}.
Since the canonical projection    $\pi : \Ghat \to H$  is an open
mapping, the image of $U$ in $H$ contains a neighborhood of the
form
\[
  C_0 = [-\epsilon /2, \epsilon /2)^d \times \{0\} \times
  \Big[-\frac{1}{2N},\frac{1}{2N}\Big)^\ell \times \{e\} \, .
\]
By construction,  $C_0$ is a fundamental domain for the lattice
\[
  \Xi = (\epsilon \bZ )^d \times \bZ ^k \times \bZ _N  ^\ell
  \times F \subseteq H .
\]
Consequently, $L^2(C_0)$ possesses an orthogonal basis consisting
of characters restricted to $C_0$ and labelled by $\Upsilon:=
\Xi^\perp $.

Since $\Xi$ is a lattice in $H$, $\Upsilon  $ is a lattice in
$\widehat{H}$. We may identify $\Upsilon $ with a subgroup of $G$
by $\Upsilon \subseteq \widehat{H} \simeq \big( \Ghat / K \big) \,
\widehat{} = K ^\perp \subseteq \widehat{\Ghat } \simeq G$.
Consequently, by fixing representatives $\widehat{k}$ from the
cosets $\widehat{K}$,  we obtain a quasi-lattice
$\Gamma=\{\widehat{k} \upsilon \}$ in $G$  with $\widehat{k}$ ranging over
$\widehat{K}$ and $\upsilon$ over $\Upsilon$.

Next, set $C = \pi \inv (C_0)$ and define for $\gamma $ in $\Gamma$
and $\omega$ in  $\widehat{G} $
\[
  \psi _\gamma (\omega ) =  \mu_{\Ghat }(C)^{-1/2}\,   \langle \omega
  , \gamma \rangle \,  \chi _{C} (\omega ) = \mu_{\Ghat
  }(C)^{-1/2} \,  \langle \omega
  , \gamma \rangle \,  \chi _{C_0} (\pi(\omega) ) \, .
\]
We now prove that the functions $\psi _\gamma$ form an orthonormal
basis for $L^2(C)$. We assume as usual that the Haar measure of a
compact subgroup $K$ is normalized to be a probability measure and
that the Haar measure of $\Ghat/K$ is normalized so that the
Weil-Bruhat formula~\cite{Rei68} $d\mu _{\Ghat } (\omega ) = d\mu
_K(k) d\mu _{\Ghat/K}(\pi (\omega ))$ holds. So we obtain that
\begin{eqnarray*}
  \mu _{\Ghat } (C) &=& \int _{\Ghat } \chi _C( \omega ) \, d\mu
  _{\Ghat }(\omega ) =  \int _{H} \int _K  \chi _C( \omega k) d\mu
  _K(k) \, d\mu
  _{H}(\pi (\omega ))\\
& =& \int _{H} \chi _{C_0}(\pi (\omega ))\, d\mu _{H}(\pi (\omega
)) = \mu _{H}(C_0)
\end{eqnarray*}
and that $\|\psi _\gamma \|_2 = 1$ for every $\gamma$ in $\Gamma$.
If $\gamma=\widehat{k}\upsilon$ and
$\gamma'=\widehat{k}'\upsilon'$ are in $\Gamma$, then using  the
Weil-Bruhat formula once more, we obtain that
\begin{eqnarray*}
\lefteqn{\int _{\Ghat } \psi _\gamma (\omega ) \overline{\psi  _{\gamma'}
(\omega )} \, d\mu _{\Ghat }(\omega )}\\
 &=& \mu _{\Ghat }(C)\inv  \int
_{H }\Big(
\int _{K } \langle \omega k, \widehat{k} \upsilon \widehat{k'}\inv
(\upsilon ')\inv \rangle \chi _C(\omega k) \, d\mu _K(k) \Big)  d\mu _{H}(\pi(\omega)) \notag \\
 & = & \delta_{\widehat{k},\widehat{k}'}
 \mu _{\Ghat }(C)\inv \int
_{H } \langle \pi (\omega ) , \upsilon (\upsilon')^{-1}
\rangle \chi _{C_0}(\pi (\omega ))\, d\mu _{H}(\pi (\omega )) = \delta
_{\gamma , \gamma' } .
\end{eqnarray*}
Here  we have used that $\langle \omega k , \upsilon (\upsilon ')\inv
\rangle $ is independent of $k$ in $K$, that $\{\langle k,
\widehat{k}\rangle \}$ is an orthonormal basis for $L^2(K)$, and that
$\{\langle \pi (\omega ),
\upsilon \rangle\}_{\upsilon \in \Upsilon}$ is an orthogonal basis
for $L^2(C_0)$.

Next we show that the linear span of $\psi _\gamma $ ($\gamma $ in
$\Gamma $) is dense in $L^2(C)$. So assume that for some $f$ in
$L^2(C)$ and  all $\gamma $ in $\Gamma $ we have
\begin{eqnarray*}
0&=&\int _{\Ghat } f(\omega ) \overline{\psi _\gamma (\omega)} \, d\mu _{\Ghat }(\omega
)\\
& =& \mu _{\Ghat }(C)^{-1/2} \int _{H} \Big( \int _K f(\omega k
) \overline{\langle \omega k , \widehat{k}\rangle}\,  d\mu _{K}(k) \Big) \,
\overline{\langle \pi (\omega), \upsilon\rangle} \chi _{C_0}(\pi
(\omega )) \, d\mu _{H}(\pi (\omega )) \, .
\end{eqnarray*}
Since $\{\langle \pi (\omega), \upsilon\rangle \}$ is an orthogonal
basis for $L^2(C_0)$, we find that
$$
 \int _K f(\omega k
) \overline{\langle \omega k , \widehat{k}\rangle } d\mu _{K}(k) = 0
$$
for almost all $\pi (\omega )$ in $\Ghat /K$ and all $\widehat{k}$ in
$\hat{K}$. We infer that $f(\omega k) = 0$ for almost all $\omega$ in
$C$ and $k$ in $K$, since $\{\langle  k , \widehat{k}\rangle \}$ is an
orthonormal basis for $L^2(K)$. Thus  the functions
 $\psi_\gamma$ form an orthonormal basis for $L^2(C)$.

To show (ii) we choose a pre-image $D$ of $\Xi$ in $\Ghat $, i.e.,
for each $\lambda $ in $\Xi$, $D\cap \pi \inv (\lambda)$ contains
exactly one element. Then $\pi (D) = \Xi$. If $dC \cap d'C \neq
\emptyset$ for $d\neq d'$ ($d,d'$ in $D$), then $\pi (d) \pi (C)
\cap\pi (d') \pi (C) = \lambda C_0 \cap \lambda ' C_0 \neq
\emptyset$ for $\lambda \neq \lambda '$ ($\lambda,\lambda'$ in
$\Xi$). Since $C_0$ is a fundamental domain for the lattice $\Xi$,
we conclude that $\lambda = \lambda '$. By choice of $D$ we also
have $d= d'$, a contradiction. Thus the translates $dC$ ($d$ in
$D$) form a partition of $\Ghat $, and (ii) is proved.
\end{proof}

\section{The Ramanathan--Steger comparison principle}
\label{sec:comparison}

The following lemma is a variation of an argument  invented by
Ramanathan and Steger~\cite{RS95}. Their decisive idea has been
investigated quite intensively in recent years.
See~\cite{GR96,CDH99,BCHL06,HK07, Kut07} for a sample of
references and \cite{Hei07} for an excellent survey. We follow the
early paper~\cite{GR96}. In what follows, $\cH$ is a separable
Hilbert space with inner product $\langle \cdot,\cdot \rangle$ and
norm $\|\cdot\|$.
\begin{lemma}\label{chap}
  Let $\Gamma$ and $\Lambda$ be uniformly discrete subsets of
  $G$. Suppose that the sequence $\{g_{\gamma} : \gamma \in \Gamma \}$
  is a Riesz sequence in $\cH $ and that there exists a sequence
  $\{h_{\lambda} : \lambda \in \Lambda \}$ so that, for fixed $\epsilon >0$
  and a compact set $K \subseteq G$,
  \begin{equation}
    \label{eq:c1}
    \mathrm{dist}\,  _{\cH} \Big( g_\gamma , \mathrm{span} \{h_\lambda :
     \lambda  \in \Lambda  \cap \gamma K\} \Big) <\epsilon
  \end{equation}
for every $\gamma\in \Gamma $. Then for every compact set
$L\subseteq G$ we have
\begin{equation}
  \label{eq:c2}
  (1-c\epsilon) \, \mathrm{card} \,  ( \Gamma \cap L ) \leq \, \mathrm{card} \,
  (\Lambda \cap L K) \, .
\end{equation}
The constant  $c>0$  depends only on $\{g_\gamma \}$. In
particular, $c=1$ if the $g_\gamma $ constitute an orthonormal
set.
\end{lemma}

\begin{proof}
  Fix a compact set $L\subseteq G $ and set
 \[
    \cH_0 = \overline{\mathrm{span} \{ g_\gamma : \gamma
      \in \Gamma \}} \, .
  \]
Then $\{ g_\gamma :\ \gamma \in \Gamma \}$ is a Riesz basis for
$\cH_0$ with dual basis $\{\tilde{g}_\gamma :\  \gamma \in \Gamma
\}\subseteq \cH_0 $, say. Since $\{\tilde{g}_\gamma \}$ is also a
Riesz basis, it is bounded, and so
\begin{equation}
  \label{eq:c4}
  c= \sup _{\gamma \in \Gamma} \|\tilde{g}_\gamma \| < \infty \, .
\end{equation}
If $\{g_\gamma :\ \gamma \in \Gamma \}$ is an orthonormal basis,
then $\tilde{g}_\gamma = g_\gamma $ and $c=1$.

Let $W_r (L) = \mathrm{span}\,\{ g_\gamma :\ \gamma \in \Gamma
\cap L \}$ and $W_f (KL) = \mathrm{span} \, \{ h_\lambda : \lambda
\in \Lambda \cap KL \}$. Let $P_{W_r}$  denote the orthogonal
projection onto $W_r(L)$ and  $Q_{W_f}$ the orthogonal projection
onto $W_f (LK)$.

Using these projections, we can recast assumption~\eqref{eq:c1} as
$\|(I-Q_{W_f})
  g_\gamma \| < \epsilon $ provided that $\gamma \in \Gamma \cap
  L$ (because in this case $\Lambda  \cap \gamma K \subseteq \Lambda
  \cap KL$). Consequently, we also have
\begin{equation}
  \label{eq:c5}
  \|(I-P_{W_r}Q_{W_f}) g_\gamma \| = \|P_{W_r}(I-Q_{W_f})P_{W_r}
  g_\gamma \| < \epsilon \qquad \text{for all }\, \gamma \in
  \Gamma \cap L
\, .\end{equation}
 The proof is done by estimating the trace of $T= P_{W_r} Q_{W_f}P_{W_r}:
 \cH_0 \to \cH_0$ in two different ways. First, since
 all eigenvalues $\nu _k$ of $T$ satisfy $0\leq \nu _k \leq 1$, we
 have
 \begin{equation}
   \label{eq:c7}
   \mathrm{tr}\, (T) \leq \mathrm{rank} \, T \leq \mathrm{dim}\, \Big(
   W_f(LK)\Big) \leq \mathrm{card} \,  \, (\Lambda \cap LK ) \, .
 \end{equation}
On the other hand, using~\eqref{eq:c4} and~\eqref{eq:c5}, we find that
\begin{eqnarray}
  \mathrm{tr}\, (T) &=& \sum _{\gamma \in \Gamma \cap L } \langle T
  g_\gamma , \tilde{g}_\gamma \rangle \notag \\
&=& \sum _{\gamma \in \Gamma \cap L } \Big(\langle
  g_\gamma ,  \tilde{g}_\gamma \rangle  -  \langle (I-T)
  g_\gamma ,  \tilde{g}_\gamma \rangle \Big) \notag \\
&\geq & \sum _{ \gamma \in \Gamma \cap L } 1 - \sum _{ \gamma
  \in \Gamma \cap L } c\epsilon \notag \\
&=& (1-c\epsilon ) \, \mathrm{card} \,  \, (\Gamma \cap L) \, .
\label{mult}
\end{eqnarray}
The claim \eqref{eq:c2} now follows from \eqref{eq:c7} and the
above.
\end{proof}

In the proof of our main theorem, we will use an orthonormal basis
with the property that $N$ functions are associated to each point
$\gamma$ in  $\Gamma$. In this case we have to count each $\gamma
$ in the final estimate ~\eqref{mult}  with multiplicity $N$. This
modification yields the following statement.

\begin{lemma}\label{chapmod}
  Let $\Gamma$ and $\Lambda $ be uniformly discrete subsets of
  $G$. Suppose that the sequence $\{g_{\gamma, j} : \gamma \in \Gamma
  , j=1, \dots , N\}$ is a Riesz sequence in $\cH $ and that there exists a sequence
  $\{h_\lambda : \lambda \in \Lambda \}$ so that, for fixed $\epsilon >0$
  and a compact set $K \subseteq G$,
  \[
    \mathrm{dist}\,  _{\cH} \Big( g_{\gamma , j} , \mathrm{span} \{h_\lambda :
     \lambda  \in \Lambda  \cap \gamma K\} \Big) <\epsilon
\]
for every $\gamma\in \Gamma $ and $j=1, \dots , N$. Then for every
compact set $L\subseteq G$ we have
\[
  (1-c\epsilon)\, N \, \mathrm{card} \,  ( \Gamma\cap L ) \leq \, \mathrm{card} \,
  (\Lambda  \cap L K) \, .
\]
The constant  $c>0$  depends only on $\{g_{\gamma,j} \}$, and
$c=1$ if $\{g_{\gamma,j} \} $ is an orthonormal set.
\end{lemma}

Our application of the Ramanathan--Steger comparison lemma will
require an estimate usually called the \emph{homogeneous approximation
property}. To state it, we introduce the following notation.
Let  $M_x$ be  the modulation operator defined by $M_x f(\omega )
= \langle \omega , x\rangle f(\omega )$ for $f\in L^2(\Ghat )$,
$x\in G, \omega \in \Ghat $.

\begin{lemma}\label{HAP}
  Let $\widehat{G}$ be compactly generated, and assume that $\{e_\lambda =
  M_\lambda g: \lambda \in \Lambda\}$, with $g$ in $L^\infty (\Omega )$,
  is a frame  for
  $L^2(\Omega )$ with dual frame $\{ h_\lambda : \lambda \in \Lambda
  \}$. Then for every $f$ in $L^2(\Omega )$ and $\epsilon >0$ there is a
  compact set $K\subseteq G$ (depending on $f$ and $\epsilon $) such
  that
  \begin{equation}
    \label{eq:c1a}
    \mathrm{dist}\,  _{\cH} \Big( M_x f , \mathrm{span} \{h_\nu
    : \nu \in \Lambda \cap x K\} \Big) <\epsilon
  \end{equation}
for every $x\in G$.
  \end{lemma}

  \begin{proof}
    The proof is identical to the proof of Lemma~2 in
    ~\cite{GR96}. Using the frame expansion of  $f\in L^2(\Omega )$,  we write
$$
M_xf = \sum _{\lambda \in \Lambda} \langle M_x f, M_\lambda
g
 \rangle  h_\lambda \, .
$$
Let $P_{x,K} $ denote the orthogonal projection from $L^2(\Omega ) $
onto $\mathrm{span}\, \{h_\lambda : \lambda\in \Lambda \cap xK\}$.
Since $\sum _{\lambda \in \Lambda \cap xK} \langle M_x f, M_\lambda
g
 \rangle  h_\lambda$ is some approximation of $f$ in $P_{x,K}L^2$,
 the square of the distance in \eqref{eq:c1a} is at most
\begin{eqnarray*}
  \|M_xf - P_{x,K}f \|_2^2 &\leq & \|\sum _{\lambda \not  \in xK} \langle
M_x f, M_\lambda  g  \rangle h_\lambda \|_2^2\\
&\leq & C \sum _{\lambda \not  \in xK} |\langle
M_x f, M_\lambda g  \rangle |^2 \\
&=& C \sum _{ \lambda \not  \in xK} |\langle
 f, M_{x^{-1}\lambda} g  \rangle |^2 \, .
\end{eqnarray*}
Set  $F(x) = \int _\Omega f(\omega ) \overline{g(\omega )}\,  \overline{\langle \omega
  ,x\rangle } d\omega = \cF \inv (f\bar{g})(x\inv)$. Then  $F\in \bom$, and  the
latter expression equals $C\sum _{\lambda \not \in xK}
|F(x^{-1}\lambda)|^2 $. If $\lambda \not \in xK$, then
$x^{-1}\lambda \not \in K$, and so we obtain as in the estimate
\eqref{eq:cd01} in the proof of Lemma~\ref{help} that
\begin{eqnarray*} \|M_xf - P_{x,K}f \|_2^2 & \le & \sum _{\lambda \not
    \in x K} \frac{1}{\mu _G(U)} \int _{x^{-1}\lambda
        U} |F^\sharp (t)|^2 d\mu _G(t)  \\
&=&  \sum _{x^{-1}\lambda \not
    \in  K} \frac{1}{\mu _G(U)} \int _{x^{-1}\lambda
        U} |F^\sharp (t)|^2 d\mu _G(t) \\
&\leq & \frac{1}{\mu _G(U)} \int _{K^c U } |F^\sharp (t)|^2 d\mu _G(t),
\end{eqnarray*}
with $U$ depending only on $\Lambda$, but not on $x\in G$. Since $F^\sharp$ is in
$L^2(G)$, we may choose $K$ so large that the expression on the
right becomes less than $\epsilon$, and this bound holds uniformly
in $x$.
  \end{proof}


%
%
%

\section{Proof of Theorem~\ref{tm:compc}}
\label{sec:proof}

The proof becomes slightly simpler if we replace $\Gamma_0$ by
\[ \Gamma_0'=(\mu_{\Ghat}(\Omega)^{1/d}\bZ)^d\times\{e\}\times D_0.
\]
This replacement can be made because it is plain that
$\Gamma_0'\preceq \mu_{\Ghat}(\Omega)\Gamma_0$ as well as $
\mu_{\Ghat}(\Omega)\Gamma_0\preceq \Gamma_0'$. Thus, by
transitivity of the relation `$\preceq$', it suffices to prove
that if the uniformly discrete set $\Lambda$ is a set of sampling
for $\bom$, then $\Gamma_0'\preceq \Lambda$, and if $\Lambda$ is a
set of interpolation for $\bom$, then $\Lambda\preceq \Gamma_0'$.

The body of the proof is an intermediate step in which we compare
$\Lambda$ with an integer multiple of one of the quasi-lattices of
Lemma~\ref{cubes}. Incidentally, this analysis applies to
$\Gamma_0'$ as well, with \[
\Omega':=[-\pi\mu_{\Ghat}(\Omega)^{1/d},\pi
\mu_{\Ghat}(\Omega)^{1/d}]^d\times\{e\}\times K.\] This
observation will enable us to eliminate the quasi-lattices. In
this part of the proof, $\Gamma_0'$ will play a ``complementary''
role to $\Lambda$; $\Gamma_0'$ is treated as a set of
interpolation for $\cB_{\Omega'}$ when $\Lambda$ is a set of
sampling for $\bom$, and vice versa.

We begin by covering $\Omega $ by an open set $\Omega_0$ such that
$\mu_{\Ghat}(\Omega _0 \setminus\Omega ) < \epsilon ^2 / 4 $. We
then take a neighborhood basis $\{ V\}$ of $e$ in $\Ghat$ and
construct the corresponding cubes $C_V$ and discrete sets
$D_V\subseteq \Ghat $ according to Lemma~\ref{cubes}. It is easy
to see that the collection $\bigcup _{V} \{ d_V C_V : d_V \in
D_V\}$ generates the Borel sets in $\Ghat $.

By taking $V$ small enough, we may choose a cube $C_0= C_V$ and a
finite number of pairwise disjoint  translates $d_j C_0$,
$d_j \in D$, $j=1, \dots , N$, such that
\[
  \Omega_* = \bigcup _{j=1}^N \Omega_j  \subseteq \Omega_0  \qquad
  \text{and } \,   \mu _{\Ghat } (\Omega  \setminus \Omega _* ) <
  \dfrac{1}{4} \epsilon^2 \mu_{\widehat{G}}(\Omega_*) \, .
\]
This is possible because the Haar measure is regular. We may even
assume that $N$ is of the form $N=2^n$ for a positive integer $n$
because the possibly discrete set of permissible values for
$\mu_{\widehat{G}}(C_0)$ is sufficiently dense. More precisely, for
arbitrary $c>1$, every interval of the form $(\delta, c\delta)$ will
contain a permissible value for $\mu_{\widehat{G}}(C_0)$ provided
that $\delta$ is sufficiently small.

By Lemma~\ref{cubes},  $L^2(d_j C) $ possesses an orthonormal
basis $\{\psi_\gamma  : \gamma \in \Gamma \}$ that is labelled by
a quasi-lattice $\Gamma $ in $G $. Consequently, $L^2(\Omega _*)$
contains an orthonormal basis of the form $\{ \psi _{\gamma , j},
\gamma \in \Gamma, j=1,\dots , N\}$ where $\psi _{\gamma , j}$ is
given explicitly by $\psi _{\gamma ,j} (\omega  ) =
\mu_{\Ghat}(C_0)^{-1/2}
  \langle \omega  , \gamma \rangle   \chi _{d_j C_0} (\pi(\omega) )$ for
  $\gamma \in \Gamma $.

We  now construct  another orthonormal basis for $L^2(\Omega _*)$  of the form
\[ \phi_{\gamma,j} (\omega  ) = \mu_{\Ghat}(\Omega_*)^{-1/2}
  \langle \omega  , \gamma \rangle  g_j (\omega )\] for
  $\gamma \in \Gamma $, where $g_j$ is a real function such that
  $|g_j|=\chi_{\Omega_*}$. We obtain $g_j$ in the following way.
Let $U= (u_{kl}), k,l = 1, \dots , N$ be a Hadamard matrix, i.e.,
$U$ has  entries $\pm 1$ and is a multiple of an orthogonal matrix.
(Such a matrix exists because $N=2^n$.) We set
\begin{equation}
 \phi _{\gamma ,j} (\omega ) = \mu _{\widehat{G}}(\Omega _*)^{-1/2} \langle
 \omega  , \gamma \rangle \sum _{k=1}^N u_{jk} \, \chi _{d_kC}(\omega
 ) \,.
\end{equation}
Then $\{ \phi _{\gamma ,j} : \gamma \in \Gamma,j=1, \dots , N\}$ is an
orthonormal basis
for $L^2(\Omega  _*)$ with $\|\phi _{\gamma ,j} \|_\infty = \mu
_{\widehat{G}}(\Omega _*)^{-1/2} $. Thus
\begin{eqnarray*}
\mathrm{dist}_{L^2} (\phi _{\gamma ,j}, L^2(\Omega )) &=& \|\phi
_{\gamma ,j}
- \phi _{\gamma ,j} \chi _\Omega \|_2 \\
& = & \|\phi _{\gamma ,j} \|_\infty \|\chi _{\Omega _*} - \chi _\Omega
\|_2 = \mu _{\widehat{G}}(\Omega _*)^{-1/2} \mu
_{\widehat{G}}(\Omega _*\Delta \Omega )^{1/2} <\dfrac{\epsilon}{2}.
\end{eqnarray*}


Let us first assume that $\Lambda$ is a set of sampling for $\bom$.
We then apply the homogeneous approximation property
(Lemma~\ref{HAP}) to the frame $e_\lambda = M_\lambda \chi _\Omega ,
\lambda \in \Lambda $, with dual frame $h_\lambda $,  and each of the
functions $g_j \chi_{\Omega}$. We then obtain a compact set $K$ such that
\[ \mathrm{dist}\, _{L^2(\widehat{G})} \Big( M_\gamma g_j \chi_{\Omega},
\mathrm{span} \{h_\lambda     \in \Lambda  \cap \gamma  K\} \Big)
< \frac{\epsilon}{2} \] 
 for $j=1, \dots , N$. Therefore,
$$
 \mathrm{dist}\, _{L^2(\widehat{G})} \Big( \phi_{\gamma,j} ,
\mathrm{span} \{h_\lambda     \in \Lambda  \cap \gamma K\} \Big) <
\epsilon \, .
$$
This is exactly the hypothesis of Lemma~\ref{chapmod}, and 
we have therefore shown that, for every compact set $L$, we have
\begin{equation}
(1-\epsilon ) N \mathrm{card}\, (\Gamma \cap L) \leq
\mathrm{card}\, (\Lambda \cap KL) \, . \label{cardinality1}
\end{equation}

If $\Lambda$ is a set of interpolation, we argue similarly. The
only difference is that now the functions $\phi_{\gamma, j}= M_\gamma g_j$ are
viewed as a frame, and the functions $e_\lambda$ constitute a Riesz
sequence. We apply again the homogeneous approximation property
and use Lemma~\ref{chap} to get
\begin{equation}\label{cardinality2} (1-c \epsilon )
\mathrm{card}\, (\Lambda \cap L) \leq N \mathrm{card}\, (\Gamma
\cap KL) \end{equation} for every compact set $L$, where $K$ is
the compact set given by Lemma~\ref{HAP}.

We have now what we need to finish the proof. To prove part (S) of
Theorem~\ref{tm:compc}, we use that $\Gamma_0'$ is a set of
interpolation for $\cB_{\Omega'}$. Hence, by \eqref{cardinality2},
there exists a compact set $K'$ such that
\begin{equation}\label{cardinality3}(1-c \epsilon ) \mathrm{card}\,
(\Gamma_0' \cap L) \leq N \mathrm{card}\, (\Gamma \cap KL)
\end{equation} holds for every compact set $L$; we may of course
adjust $\epsilon$ and the approximation of $\Omega$ so that the
$\Gamma$ also suits the approximation of $\Omega'$. If $\Lambda$ is
a set of sampling, then combining \eqref{cardinality1} with
\eqref{cardinality3}, we obtain that
$$
(1-c \epsilon ) \mathrm{card}\,
(\Gamma_0' \cap L) \leq \frac{1}{1-\epsilon} \mathrm{card}\, (\Lambda
\cap K^2L) \, ,
$$
from which the desired relation $\Gamma_0' \preceq \Lambda$
follows.

Reversing the roles of $\Lambda$ and
$\Gamma_0'$, we obtain similarly $\Lambda\preceq \Gamma_0'$ when
$\Lambda$ is a set of interpolation for $\bom$.

\section{Properties of uniform densities}\label{sec:dense}

We return to some basic questions about uniform densities that
were raised in Section~\ref{sec:lcagroups}.

\begin{lemma}
  \label{densities}
For every uniformly discrete subset $\Lambda$ of an LCA group $G$,
we have  $\mathcal{D}^-(\Lambda)\le \mathcal{D}^+(\Lambda
)<\infty$.
\end{lemma}

\begin{proof}
It is sufficient to prove that both $\cD^+(\Lambda)<\infty$ and
$\cD^-(\Lambda)<\infty$. Indeed, if $\cD^-(\Lambda)>
\cD^+(\Lambda)$, then $\Lambda\preceq \delta \Lambda$ for some
$\delta<1$. By the transitivity of the relation `$\preceq$', this
can only happen if $\cD^-(\Lambda)=0$ or $\cD^{-}(\Lambda)=\infty$.

We first prove that $\cD^+(\Lambda)<\infty$. We need to show that
there exists a positive number $\alpha$ such that $\Lambda \preceq
\alpha \Gamma_0$. Let $L$ be a compact subset of $\Lambda$. Since
$\Lambda$ is uniformly discrete, there is a uniform bound, say
$M$, on the number of points from $L\cap \Lambda$ to be found in
each set $\gamma K$, where $K:=[-1/2,1/2]^d \times \bT^m \times\{e\}$
and $\gamma$ is an element in $\Gamma_0$. Therefore,
\[ \mathrm{card}\, (\Lambda \cap L )
  \leq  M \mathrm{card}\, (\Gamma_0 \cap KL), \]
  and so $\Lambda \preceq
M \Gamma_0$

We next prove that $\cD^-(\Lambda)<\infty$. Let us assume that we
have $\alpha \Gamma_0 \preceq \Lambda$ for some $\alpha$. Then for
every positive $\epsilon$ there exists a compact set $K$ such that
\begin{equation}
  \label{eq:la1}
  (1-\epsilon) \, \alpha \, \mathrm{card}\, (\Gamma_0 \cap L )
  \leq  \mathrm{card}\, (\Lambda \cap KL)
\end{equation}
  for every compact set $L$.
  We may assume that $K=B\times \bT^m\times F$, where $B$ is a ball in $\rd$ centered
  at the origin and $F$ is a finite subset of $D_0$ such $F^{-1}=F$.
  Then $\bigcup_{n=1}^\infty F^n$ is a finitely generated subgroup of
  $D_0$, which has the structure $\Z^l\times E$ with $E$ a
  finite group. (See \cite[p. 451]{HR63}.) To simplify the argument, we may assume that
  $F$ is just $B'\times E$, with $B'$ a ball in $\Z^l$ centered at the origin.
We choose $L=K^n$ and note that for sufficiently large $n$ we have
\begin{equation}
  \label{eq:la2}
   \mathrm{card}\, (\Gamma_0 \cap L )\, \ge (1-\epsilon)\, \mu_G(L).
\end{equation}
On the other hand, if $U\subseteq K$ is an open set  such that the sets
$\lambda U$ ($\lambda$ in $\Lambda$) are pairwise disjoint, we
obtain
\begin{equation}
  \label{eq:la3}
   \mathrm{card}\, (\Lambda  \cap KL) \le \mu_G(U)^{-1}
\mu_G(K^{n+2}) \le (1+\epsilon) \,\mu_G(U)^{-1} \mu_G(L)
\end{equation}
whenever $n$ is sufficiently large. Combining \eqref{eq:la1} --
\eqref{eq:la3}, we obtain that for $\epsilon >0$
$$
\alpha \leq \frac{1+\epsilon}{(1-\epsilon)^2}\,  \mu _G(U)\inv \, ,
$$
and thus $\cD^{-}(\Lambda)\le \mu_G(U)^{-1}$.
\end{proof}

The relation `$\preceq$' may be viewed as a relation between
discrete measures. Since the canonical lattice $\Gamma_0$ has a
highly regular distribution, it should come as no surprise that we
may replace the discrete measure associated with $\Gamma_0$ by the
Haar measure $\mu_G$. Interpreting a uniformly discrete set as a
sum of point masses located at the points $\lambda$ of the set, we
may generalize the relation `$\preceq$' to arbitrary nonnegative
measures on $G$. Thus, if $\nu$ and $\tau$ are two such measures
on $G$, we write $\nu\preceq\tau$ if for every $\epsilon>0$ there
exists a compact set $K$ in $G$ such that \[ (1-\epsilon)\nu(L)\le
\tau(LK)\] for every compact set $L$ in $G$. If we set  again \[
K=[-1/2,1/2]^d\times \bT^m\times \{e\}, \] then it is immediate
that \[ \mu_G(L)\le \mathrm{card}\, (\Gamma_0 \cap KL ) \ \ \
\text{and} \ \ \  \mathrm{card}\, (\Gamma_0 \cap L ) \le
\mu_G(KL)\] for every compact set $L$. This implies that
$\mu_G\preceq\Gamma_0$ and $\Gamma_0\preceq \mu_G$,  so that
Theorem~\ref{tm:compc} can be restated in the following form.

\begin{Th1primeprime}
Suppose $\Lambda$ is a uniformly discrete subset of the LCA group
$G$ and $\Omega$ is a relatively compact subset of $\widehat{G}$.
\begin{itemize}
\item[(S)] If $\Lambda $ is a set of sampling for $\bom $, then
$ \mu_{\Ghat}(\Omega)\, \mu_G \preceq \,
 \Lambda$. \item[(I)] If $\Lambda $ is a set of interpolation
for $\bom $, then  $  \Lambda\,  \preceq \,
 \mu_{\Ghat}(\Omega)\, \mu_G $.
\end{itemize}
\end{Th1primeprime}

We finally show that, in $\rd $, our uniform densities coincide with
the classical Beurling densities. In $\rd $ we use the standard
additive notation $x+y$ and $K+L$ instead of the multiplicative
notation on arbitrary LCA groups, and we write $|U|$ for the
Lebesgue (Haar) measure of $U\subseteq \rd $.

\begin{lemma}
  \label{equiv}
If $G=\rd $, then $\mathcal{D}^-(\Lambda ) =
\mathcal{D}^-_B(\Lambda )$ and  $\mathcal{D}^+(\Lambda ) =
\mathcal{D}^+_B(\Lambda )$ for every uniformly discrete set
$\Lambda$.
\end{lemma}

\begin{proof}
Let $\Lambda$ be a uniformly discrete subset of $\rd$. Then, for
every $\epsilon>0$, there exists a compact set $K=Q_R(0)=
[-R/2,R/2]^d$ such that
\[
   (1-\epsilon )  \cD
^-(\Lambda ) \mathrm{card} \, (\Z^d \cap L) \leq  \mathrm{card}
   \, \big(\Lambda \cap (L+Q_R(0))\big) \,
\]
for every compact set $L$. Specializing to cubes  $L=Q_h(y), y\in
\rd ,$ we get that
$$
   (1-\epsilon ) \cD
^-(\Lambda )  \inf _{y\in \rd } \frac{\mathrm{card} \, (\Z ^d
     \cap Q_h(y))}{(h+R)^d} \leq \inf _{y\in \rd } \frac{\mathrm{card}
   \, (\Lambda \cap Q_{h+R}(y))}{(h+R)^d} \, .
$$
Taking the limit $h\to \infty $, we obtain $(1-\epsilon)\cD
^-(\Lambda )\le \cD _B^-(\Lambda )$, and so $\cD ^-(\Lambda )\le
\cD _B^-(\Lambda )$ since the inequality holds for every positive
$\epsilon$.

Conversely, we may for any given $\epsilon
>0$ find  $h_0>0$ such that
\begin{equation}
  \label{eq:cd4}
  \frac{\mathrm{card}\, (\Lambda \cap Q_h(y))}{h^d} \geq (1-\epsilon
  )\, \cD ^-_B
(\Lambda )
\end{equation}
for every point $y$ in $\rd$ and $h>h_0$. Now partition $\rd $
into cubes $Q_h(hk), k\in \zd, $ whose interiors are disjoint.
Given a compact set $L\subseteq \rd $, there exist finitely many
$k_j \in \zd, j=1,\dots  J$, such that
\[ L \subset \bigcup_{j=1}^J Q_h(hk_j) \subset L+Q_{2h}(0)\, . \]
Then by \eqref{eq:cd4}
\[
\mathrm{card}\, \Big(\Lambda \cap (L+Q_{2h}(0))\Big)\, \geq \, \sum _{j=1}^J
\mathrm{card}\, (\Lambda \cap Q_h(hk_j)) \\
\, \geq \,  h^d (1-\epsilon ) \,J\,\cD ^-_B (\Lambda ).
\]
Since $(h+1)^d\ge \mathrm{card}\, ( \Z^d\cap Q_h(hk_j))$, it follows
that
\[ \mathrm{card}\, \big(\Lambda \cap (L+Q_{2h}(0))\big) \, \ge \,
(1-\epsilon )\,(1+1/h)^{-d}\,\cD^-_B(\Lambda)\, \mathrm{card}\,
(\Z^d \cap L)\, .\] We may choose $\epsilon$ arbitrarily
small and $h$ arbitrarily large, and hence $\cD ^-(\Lambda ) \geq
\cD ^-_B (\Lambda )$.

The identity  $\mathcal{D}^+(\Lambda ) =
\mathcal{D}^+_B(\Lambda )$ is proved similarly.
\end{proof}

\section{Arbitrary LCA Groups} \label{sec:red}
So far we have assumed that $\widehat{G}$ is compactly generated. 
This is not a serious restriction, as shown by the following
lemma. (See also~\cite{fg92}.)

\begin{lemma}\label{wlog}
Assume that $\Omega \subseteq \widehat{G}$ is relatively compact and
let $H$ be the open subgroup generated by $\Omega \subseteq
\widehat{G}$. Then $H$ is compactly generated and
  there exists a compact subgroup $K \subseteq G$ such that every $f
  \in \bom $ is $K$-periodic.

 Furthermore, the quotient $G/K$ factors
  as $G/K \simeq \rd \times \bT ^k\times D_0 $ for some countable
  discrete abelian group $D_0$ and $(G/K)\, \widehat{}\, = H$, where $H$
  is the  open subgroup of $\Ghat $ that is generated by  the spectrum $\Omega
  $.
\end{lemma}

\begin{proof}
   Choose an open, relatively compact neighborhood $V$  of the
   spectrum $\Omega \subseteq \Ghat$, and let $H $ be the
   \emph{open} subgroup of $\Ghat$ that is generated by $V$. Then
   $\Ghat / H $ is discrete, and thus the group $ \big( \Ghat / H
   \big)\, \widehat{}\, $ is compact. We claim that $K:=H^\perp  $ is
   the subgroup we are looking for. Let $f\in \bom $, $x\in G, k\in
   K$, then  by the inversion formula
   \begin{eqnarray*}
f(xk) &=& \int _{\Omega} \widehat{f} (\omega ) \omega (xk) \, d\mu
_{\Ghat}(\omega ) \\
&=& \int _{\Omega} \widehat{f} (\omega ) \omega (x) \omega (k) \,
d\mu
_{\Ghat}(\omega )\\
&=& \int _{\Omega} \widehat{f} (\omega ) \omega (x) \, d\mu
_{\Ghat}(\omega ) =  f(x)
   \end{eqnarray*}
since $k\in H^\perp $ and $ \Omega \subseteq H$.

Since $H$ is compactly generated, $H$ is isomorphic to a group $H
\simeq \rd \times \bZ ^k \times L$ for some compact group $L$ by the
structure theorem for LCA groups~\cite[Thm.~9.8]{HR63}.
Consequently,
$$
\widehat{H}\, \simeq \widehat{\Ghat} / H^\perp \simeq G/K \simeq
\rd \times \bT ^k \times D_0 \, ,
$$
where $D_0 = \widehat{L}$ is a  discrete group.
\end{proof}

Consequently, every bandlimited function $f\in \bom $ lives on a
quotient $G/K$ and  may be identified with a function $\tilde{f}
\in L^2(G/K)$. \medskip

\noindent \textbf{Example}. Let $\bQ _p$ be the group of $p$-adic
numbers~\cite{HR63} with dual
  group isomorphic to $\bQ _p$.  The
  $p$-adic numbers possess a ``quasi-metric'' $|\cdot |_p$ such that
  $|x+y|_p \leq \max (|x|_p,|y|_p)$ for all $x,y \in \bQ
  _p$. Moreover,  for
  each $n\in \bZ $,
$K_n := \{ x\in \bQ _p : |x|_p \leq n\}$  is a compact-open subgroup
of $\bQ _p$. As a consequence, every relatively compact set $\Omega
\subseteq \bQ _p$ generates a compact group $H$ contained in some
$K_n$. In particular, $\bQ _p$ does not contain any  lattice.

It seems that our main theorem does not say anything about
sampling in $p$-adic groups. However,  Lemma~\ref{wlog} says that
we may assume without loss of generality that $\Ghat $ is one of
the $K_n$'s where $K_n$ contains the group $H$ generated by the
spectrum $\Omega $. Furthermore,  all functions in $\bom $ are
$H^\perp $-periodic and thus live on the discrete group $\bQ
_p/H^\perp $.  Thus we may apply Theorem~\ref{tm:compc} to the
pair $G=\bQ_p/ H^\perp $ and $H\subseteq K_n$.

\section{Closing remarks}\label{sec:fin}

(\textbf{1}) In his paper \cite{Lan67}, Landau made a slightly
weaker assumption on $\Omega$ when considering sets of sampling.
Instead of taking $\Omega$ to be relatively compact, he assumed
that $\Omega$ had positive measure. It is clear that we may
similarly take $\Omega$ to have positive Haar measure in part (S)
of Theorem~\ref{tm:compc} because such $\Omega$ can be
approximated by compact sets contained in $\Omega$. Note that this
relaxation cannot be made in part (I) of Theorem~\ref{tm:compc}.

\medskip

(\textbf{2}) Landau used his results in \cite{Lan67} to prove a
conjecture of A. Beurling concerning the lower uniform density of
sets in $\rd$ for which so-called balayage is possible. We do not
wish to go into detail about Beurling's problem, but we would like
to point out that, using our notion of density, we may extend
Landau's result concerning balayage. The restriction we have to
make is that the group $G$ be of the form $G=\rd\times\Z^m\times
K_0$ with $d\ge 1$. Theorem~5 in \cite{Lan67} extends from the
setting of $\rd$ to such groups, under the same regularity
conditions on the spectrum. The details needed to carry out this
extension can be found in \cite[pp. 341--350]{Be89} and in
Landau's paper \cite{Lan67}.

\medskip

(\textbf{3}) In his thesis \cite{Mar08}, Marzo proved that for
every relatively compact set $\Omega$ in $\rd$ we can find sets of
sampling and sets of interpolation for $\bom$ of Beurling
densities arbitrarily close to those given by Landau's theorem. It
would be interesting to know if, similarly, our density conditions
are optimal for every relatively compact set in a general LCA
group.

\medskip

(\textbf{4}) In section~\ref{sec:lcagroups}, we excluded the case
of compact groups. Our result is certainly of no interest for
compact groups, but for such groups one can state closely related
and nontrivial problems. An example is the recent work of J.
Ortega-Cerd\`{a} and J. Saludes  on Marcinkiewicz-Zygmund
inequalities \cite{OS07}. Their work deals with the group $G=\bT$
and the asymptotic behavior of sets of sampling and interpolation
when the size of the spectrum grows and we require uniform bounds
on the norms. Another, probably much more difficult problem, is to
describe similarly asymptotic density conditions when $G=\bT^m$
and both the spectrum and $m$ grow.


 \bibliographystyle{abbrv}
 \bibliography{general,new}

\end{document}

%% file: macros.tex
\newcommand{\fif}{if and only if}

\newcommand{\onb}{orthonormal basis}

\newtheorem{tm}{Theorem}
\newtheorem*{Th1prime}{Theorem 1'}
\newtheorem*{Th1primeprime}{Theorem 1''}
\newtheorem{lemma}[tm]{Lemma}
\newtheorem{prop}[tm]{Proposition}

\newcommand{\rem}{\noindent\textsl{REMARK:}}

 \theoremstyle{definition}
 
\newtheorem{remk}{Remark}

\newtheorem{exmps}{Examples}

\newcommand{\beqa}{\begin{eqnarray*}}
\newcommand{\eeqa}{\end{eqnarray*}}

\DeclareMathOperator*{\supp}{supp}

\newcommand{\field}[1]{\mathbb{#1}}
\newcommand{\bR}{\field{R}}        
\newcommand{\bN}{\field{N}}        
\newcommand{\bZ}{\field{Z}}        
\newcommand{\bQ}{\field{Q}}        
\newcommand{\bT}{\field{T}}        %
\def\La{\Lambda}

 \def\cF{\mathcal{F}}              
 
 \def\cD{\mathcal{D}}
 \def\cH{\mathcal{H}}
 \def\cB{\mathcal{B}}

 \def\cA{\mathcal{A}}

\def\rd{\bR^d}

\def\lrd{L^2(\rd)}

\def\zd{\bZ^d}

\def\<{\left<}
\def\>{\right>}

\def\inv{^{-1}}

\def\mv1{M_v^1}

